\documentclass[11pt]{article}

\usepackage{graphicx}
\usepackage{amsmath,amsthm,amssymb,enumerate}
\newtheorem{theorem}{\bf Theorem}[section]

\newtheorem{definition}{\bf Definition}[section]
\newtheorem{problem}{\bf Problem}[section]

\newtheorem{proposition}{\bf Proposition}[section]
\newtheorem{remark}{\bf Remark}[section]
\newtheorem{lemma}{\bf Lemma}[section]

\makeatletter
 \renewcommand{\theequation}{%
   \thesection.\arabic{equation}}
  \@addtoreset{equation}{section}
\makeatother
\makeatletter
\@addtoreset{equation}{section}
\renewcommand{\theequation}{\thesection.\@arabic\c@equation}
\makeatother

\begin{document}

\title{On extracting the positions of multiple unknown cracks that occur on the junction line of two elastic plates}
\author{Masaru IKEHATA\footnote{
Laboratory of Mathematics, Graduate School of Advanced Science and Engineering, Hiroshima University, Higashihiroshima 739-8527, JAPAN.
Emeritus Professor at Gunma University, Maebashi 371-8510, JAPAN.
E-mail:ikehataprobe@gmail.com
},
and Hiromichi ITOU\footnote{Department of Mathematics,Tokyo University of Science,Tokyo 162-8601,JAPAN.
E-mail:h-itou@rs.tus.ac.jp
}
}
\date{}
\maketitle

\begin{abstract}
This paper is concerned with the reconstruction issue of an inverse crack problem in a two-dimensional bounded domain which may have a possible application to the nondestructive evaluation
of materials.
It is assumed that the domain consists of two elastic plates welded together and has some unknown cracks on the junction line and the governing equation of in-plane displacement is Navier's equation.
The problem is to extract information about the location of cracks from the observation data which is a single set of a loading surface traction and the resulted in-plane displacement field on the boundary of the domain.
It is shown that the enclosure method combined with the Kelvin transform yields explicit extraction formulae of such information from the observation data.
\end{abstract}

\noindent
MSC 2010: 35R30, 74B05.

\noindent KEY WORDS: inverse crack problem, enclosure method, linearized elasticity.

\section{Introduction}
We deal with a reconstruction problem for multiple cracks located on a single line in a linearized elastic plate.
This happens when using spot welding to join two elastic plates to form a plate. The out of welded part on the faying surface can be considered a set of cracks, e.g. \cite{J2002}.
Developing mathematical methods for estimating such parts from the data observed at the plate boundaries has the potential to be applied to non-destructive evaluation of materials.

As one of mathematically exact methods, in \cite{HIIS, IIS} we have already developed a method which employs the enclosure method \cite{I1999, I2003} combined with the idea of using the Kelvin transform back to \cite{I2005} in the case of electric conductive plate.
The aim of this paper is to extend the results in \cite{HIIS, IIS} to the linearized elastic plate case.
The main parts consist of two theorems for extracting information about the location of the unknown cracks.
The one is an extension of Theorem 1 in \cite{IIS}.
The second is an extension of Theorem 2.1 in \cite{HIIS},
which already announced in a conference report \cite{HI2020} briefly. However, we will supplement the part that was not described in \cite{HI2020}, which is an application of the idea of taking the logarithmic derivative of
the indicator function in the enclosure method developed in \cite{I2011}.  See also \cite{I2021} in which this idea
has been applied to an inverse source problem.  

It should be noted that, in \cite{HS} a novel multi-modality fusing electrical and elasticity imaging is proposed and the possibility to stabilize the inversion process is suggested by complementing information obtained from both modalities each other.
We share the idea that better information can be obtained by utilizing the data obtained by measuring different physical quantities of one material in different methods.

The rest of the paper is organized as follows.
In Section 2, a formulation of the corresponding forward problem is given and then a crack detection problem which we consider in this paper is described.
In Section 3, we introduce mathematical tools in the enclosure method for solving our problem and state our main results, that is, Theorem \ref{th0} and \ref{th1}. 
Section 4 is devoted to proving the proof of Theorem \ref{th0} and \ref{th1}.
In Section 5, we describe the proof of a proposition and a lemma which play crucial role in that of those theorems.
Finally in Section 6, concluding remarks are given.

\section{Formulation}

First we describe the geometry of a material made by joining two elastic plates.
Choose the two-dimensional Cartesian coordinates $(x_1,x_2)$ in such a way that the region $\Omega$ where the material occupied
takes the form $\Omega=]0,\,a[\times\,]0,b[$ with positive numbers $a$ and $b$.
We suppose all the multiple cracks lie on the segment $[0,\,a]\times\{c\}$ with a fixed $c\in\,]0,\,b[$.  
This means that $\Omega$ is obtained by joining two elastic plates $\Omega^{+}=\Omega\cap\,\{\mbox{\boldmath $x$}=(x_1,x_2)\in{\mathbb R}^2\,\vert\,x_{2}>c\}$
and $\Omega^{-}=\Omega\cap\{\mbox{\boldmath $x$}=(x_1,x_2)\in{\mathbb R}^2\,\vert\,x_{2}<c\}$ and the junction line which is the one-dimensional version of the faying surface,
is given by $x_2=c$ and $0\le x_1\le a$.
Denote by $\Sigma$ the set of all cracks in $\Omega$ and $x_{1}$-components of crack tips by $c_{0}<c_{1}<\cdots<c_{2m+1}$ with $m\ge 1$.  
Thus the $\Sigma$ takes the form
\[
\Sigma =\bigcup_{j=0}^{m}\,[c_{2j},\,c_{2j+1}]\times\{c\}.
\]

In this paper, same as \cite{HIIS, IIS} we assume that $c_0=0$ and $c_{2m+1}=a$.
This means that the leftmost and rightmost cracks are exposed on the surface.
Note that other cases can be treated without any essential change.

\begin{figure}[htbp]
\centering
\includegraphics[width=6.5cm]{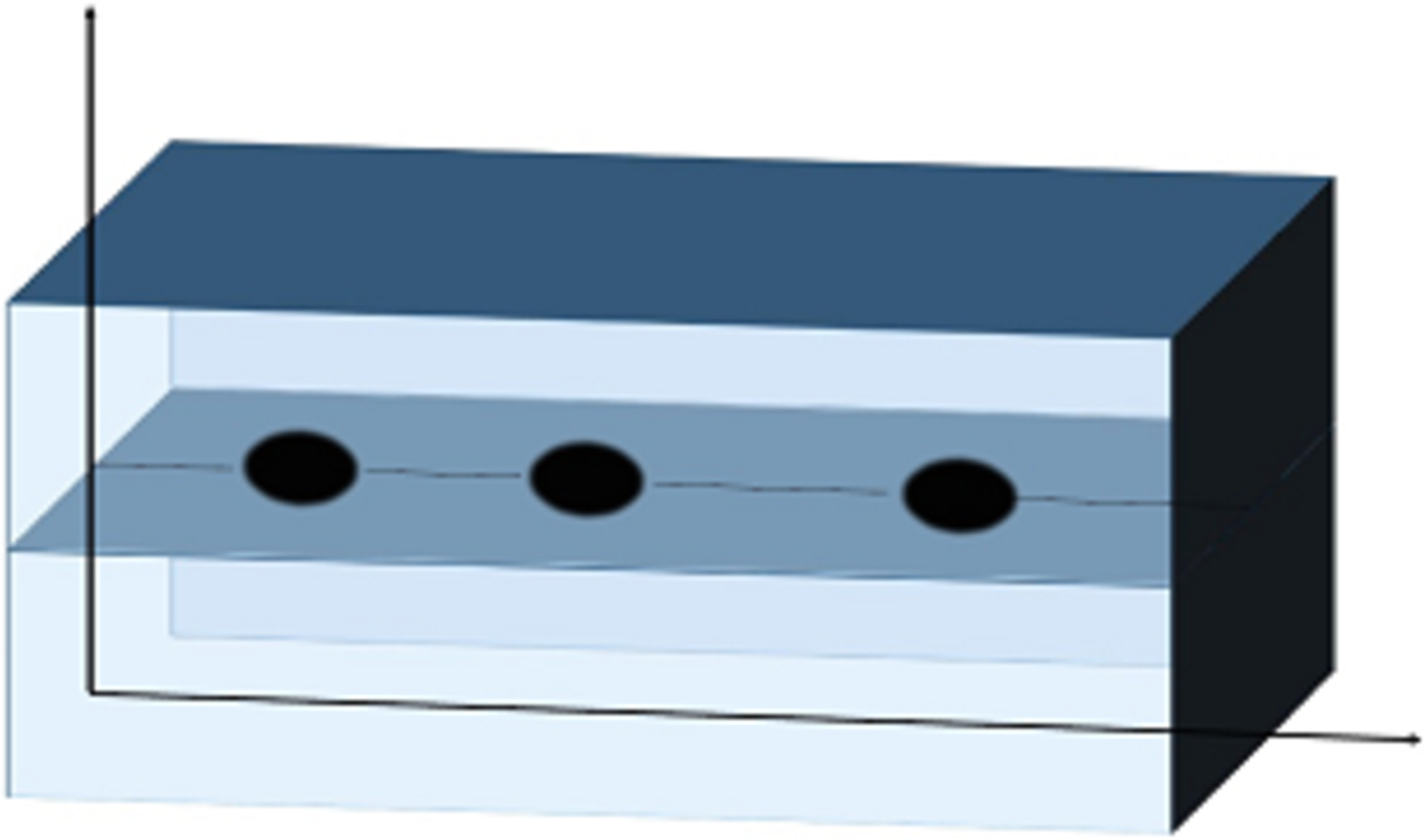}
$\;$
\includegraphics[width=6.5cm]{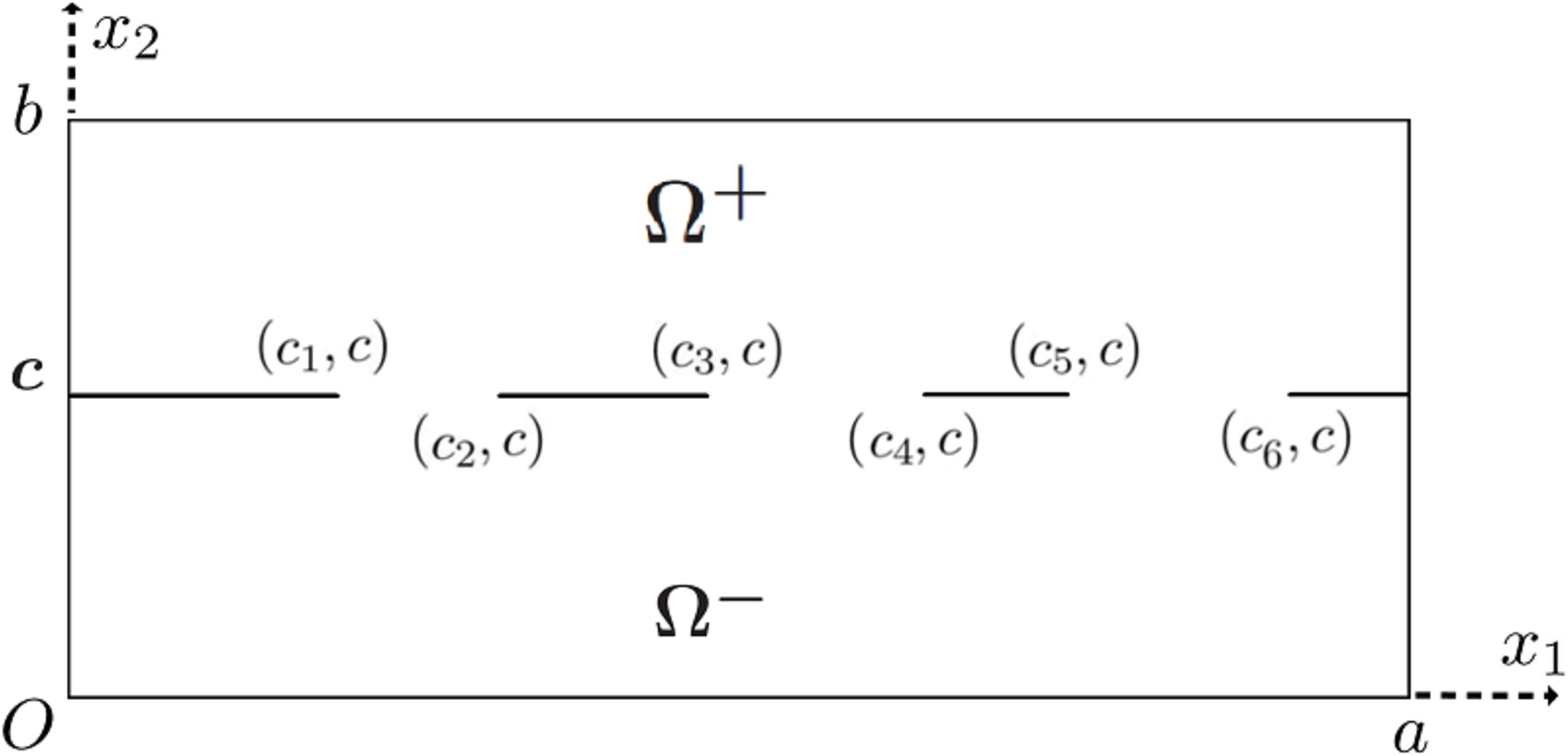}
\caption{(left): an illustration of a joined two elastic plates by spot welding whose parts indicated by filled black ellipses.
(right): an example of the domain $ \Omega $}
\label{fig1}
\end{figure}

Needless to say, our original problem should be formulated in three dimensions such as Figure \ref{fig1} (left) and thus the junction line becomes the faying surface.
However, as the first attempt, we here consider its two-dimensional version such as the cross section illustrated in Figure \ref{fig1} (right).

Next we introduce the linearized elasticity equation and the corresponding boundary value problem.
Let $ \mbox{\boldmath $u$}=(u_i)_{i=1,2} $, $ \mbox{\boldmath $\varepsilon$} =(\varepsilon_{ij})_{i,j=1,2} $ and $ \mbox{\boldmath $\sigma$} =(\sigma_{ij})_{i,j=1,2} $ be the displacement vector, the linearized strain tensor and the stress tensor, respectively.
The relation between $ \mbox{\boldmath $\varepsilon$} $ and $ \mbox{\boldmath $u$} $ is given by
\begin{eqnarray}
\mbox{\boldmath $\varepsilon$}(\mbox{\boldmath $u$}) = \frac{1}{2} \left( \nabla \mbox{\boldmath $u$} +(\nabla \mbox{\boldmath $u$})^{\rm T}\right) \label{eq:1-3}
\end{eqnarray}
where the superscript ${\rm T}$ denotes matrix transposition.
For the linearized elasticity which is a homogeneous isotropic body in the state of a plane strain, the constitutive law, so-called Hooke's law, is described as follows
\begin{eqnarray}
\mbox{\boldmath $\sigma$}=\lambda {\rm tr} (\mbox{\boldmath $\varepsilon$}) \mbox{\boldmath $I$}+2\mu \mbox{\boldmath $\varepsilon$} \label{eq:1-1}
\end{eqnarray}
where $\mbox{\boldmath $I$}$ is the identity tensor, $ \lambda $ and $ \mu $ are $ Lam\acute{e} $ constants satisfying $\mu >0$ and $\lambda +\mu >0$.
From the law of conservation of momentum the static equilibrium equation in the absence of body forces becomes
\begin{eqnarray}
\nabla\cdot \mbox{\boldmath $\sigma$} =\mbox{\boldmath $0$}. \label{eq:1-2}
\end{eqnarray}
Substituting \eqref{eq:1-1} and \eqref{eq:1-3} into \eqref{eq:1-2}, we arrive at the governing equations for $ \mbox{\boldmath $u$}=(u_1,u_2) $
\begin{eqnarray}
\mu\triangle\mbox{\boldmath $u$} + (\lambda +\mu) \nabla (\nabla\cdot\mbox{\boldmath $u$})=\mbox{\boldmath $0$}. \label{eq:1-4}
\end{eqnarray}

For given $ \mbox{\boldmath $g$}\in L^{2}(\partial\Omega) $ which is the surface force acting on $\partial\Omega$ we consider the following boundary value problem $(*)$
\[
(*)\left\{ \begin{array}{l}
\vspace{0.3cm}
\mu\triangle\mbox{\boldmath $u$} + (\lambda +\mu) \nabla (\nabla\cdot\mbox{\boldmath $u$})=\mbox{\boldmath $0$} \quad {\rm in} \quad \Omega\setminus\Sigma,\\
\vspace{0.3cm}
\mbox{\boldmath $\sigma$}^{+}\mbox{\boldmath $\nu$} =\mbox{\boldmath $\sigma$}^{-}\mbox{\boldmath $\nu$} =\mbox{\boldmath $0$} \quad {\rm on} \quad \Sigma^{\pm},\\
\mbox{\boldmath $\sigma$}\mbox{\boldmath $\nu$}=\mbox{\boldmath $g$}\quad {\rm on} \quad \partial\Omega .
\end{array}\right.
\]
Here on the crack $\Sigma$ the free traction condition is imposed
where the upper and lower sides of the stress vector are denoted by  $\mbox{\boldmath $\sigma$}^{+}\mbox{\boldmath $\nu$} $ and $ \mbox{\boldmath $\sigma$}^{-}\mbox{\boldmath $\nu$} $ with a fixed normal vector $ \mbox{\boldmath $\nu$}=(0,1) $ on $\Sigma$.
On the boundary of $\Omega$ we assume the standard Neumann type boundary condition with the unit outward normal $ \mbox{\boldmath $\nu$} $.
The governing equation $(*)$ means that both $\Omega^+$ and $\Omega^-$ consist of the same isotropic homogeneous elastic material.

Let ${\cal R}$ be the space of rigid displacements described as 
\begin{eqnarray*}
{\cal R}=\{ \mbox{\boldmath $v$}\; |\; \mbox{\boldmath $v$}( \mbox{\boldmath $x$}) =(k_1 +k_0 x_2 ,k_2 -k_0 x_1),\; k_0 ,k_1, k_2\in {\mathbb R}\}.
\end{eqnarray*}
We employ the variational formulation of $(*)$ and define the weak solution as follows.
For given $ \mbox{\boldmath $g$}\in L^{2}(\partial\Omega) $ satisfying
\begin{eqnarray}
\forall \mbox{\boldmath $\rho$}(\mbox{\boldmath $x$}) \in {\cal R},\quad \int_{\partial\Omega}\mbox{\boldmath $g$}\cdot \mbox{\boldmath $\rho$}\; {\rm d}s_{\mbox{\boldmath $x$}}=0, \label{eq:1-6}
\end{eqnarray}
we call $ \mbox{\boldmath $u$}\in H^{1}(\Omega\setminus\overline{\Sigma})/ {\cal R} $ a weak solution of the problem $(*)$ if for arbitrary $ \mbox{\boldmath $\varphi$}\in H^{1}(\Omega\setminus\overline{\Sigma})/ {\cal R} $ it holds
\begin{eqnarray}
\int_{\Omega\setminus\overline{\Sigma}} \mbox{\boldmath $\sigma$}(\mbox{\boldmath $u$}): \mbox{\boldmath $\varepsilon$}(\mbox{\boldmath $\varphi$}) \;{\rm d}\mbox{\boldmath $x$} = \int_{\partial\Omega}\mbox{\boldmath $g$}\cdot \mbox{\boldmath $\varphi$}\; {\rm d}s_{\mbox{\boldmath $x$}}, \label{eq:2-1}
\end{eqnarray}
where the double dot in $ \mbox{\boldmath $\sigma$}(\mbox{\boldmath $u$}): \mbox{\boldmath $\varepsilon$}(\mbox{\boldmath $\varphi$})= \Sigma_{i,j=1,2} \sigma_{ij}(\mbox{\boldmath $u$})\varepsilon_{ij}(\mbox{\boldmath $\varphi$})$ implies the scalar product of matrices.
It is well-known that there exists a unique weak solution of the problem ($ * $).

In this paper, we consider the following crack detection problem.
\begin{problem}
Apply the surface force $\mbox{\boldmath $g$}\not=\mbox{\boldmath $0$}$
satisfying \eqref{eq:1-6} and measure the corresponding displacement field $\mbox{\boldmath $u$}$ on $\partial\Omega$. 
Extract information about the exact location of $\Sigma$ from the singe set of data $\mbox{\boldmath $g$}$ and $\mbox{\boldmath $u$}$ on $\partial\Omega$. 
\end{problem}

\section{Statement of the main results}

First we introduce a special solution of \eqref{eq:1-4} in a neighbourhood of $\overline{\Omega}$.
Given an arbitrary point $\mbox{\boldmath $x$}\in{\mathbb R}^2\setminus\overline{\Omega}$ define
\footnote{One may choose another one 
$\mbox{\boldmath $v$}_{\tau}(\mbox{\boldmath $y$};\mbox{\boldmath $x$})=\nabla_{\mbox{\boldmath $y$}}\left(v_{\tau}(\mbox{\boldmath $y$};\mbox{\boldmath $x$})\,\right)$.  
This vector valued function also satisfies equation \eqref{eq:1-4} and the divergence free property.
However, the computation of the asymptotic behaviour of the indicator function defined later shall become complicated.
}
\[
\mbox{\boldmath $v$}_{\tau}(\mbox{\boldmath $y$};\mbox{\boldmath $x$}) := (\mbox{\boldmath $e$}_1+{\rm i}\mbox{\boldmath $e$}_2)v_{\tau}(\mbox{\boldmath $y$};\mbox{\boldmath $x$}), \quad
\mbox{\boldmath $y$}\in{\mathbb R}^2\setminus\{\mbox{\boldmath $x$}\},
\]
where the function $v_{\tau}(\mbox{\boldmath $y$};\mbox{\boldmath $x$})$ is given by
$$\begin{array}{ll}
\displaystyle
v_{\tau}(\mbox{\boldmath $y$};\mbox{\boldmath $x$})
&
\displaystyle
:=\left. e^{-\tau \mbox{\boldmath $z$}\cdot(\mbox{\boldmath $e$}_2+{\rm i}\mbox{\boldmath $e$}_1)}\right|_{\mbox{\boldmath $z$}=\frac{\mbox{\boldmath $y$}-\mbox{\boldmath $x$}}{\vert \mbox{\boldmath $y$}-\mbox{\boldmath $x$}\vert^2}}
\\
\\
\displaystyle
&
\displaystyle
=\exp\left\{-\frac{{\rm i}\tau}{(y_1-x_1)+{\rm i}(y_2-x_2)}\right\},
\end{array}
$$
$\mbox{\boldmath $e$}_{1}:=(1,0)$, $\mbox{\boldmath $e$}_2:=(0,1)$, ${\rm i}:=\sqrt{-1}$ and $\tau$ is a positive (large) parameter.
The $v_{\tau}$ is nothing but the Kelvin transformation of the complex geometrical optics solution of the Laplace equation which has been used in \cite{HIIS, IIS}.
Since $v_{\tau}$ is a {\it holomorphic} function of the complex variable $y_1+{\rm i}y_2\,(\not=x_1+{\rm i}x_2)$, 
the vector valued function $\mbox{\boldmath $v$}_{\tau}$ automatically satisfies the divergence free property and clearly equation \eqref{eq:1-4} in the domain ${\mathbb R}^2\setminus\{\mbox{\boldmath $x$}\}$.
Besides, for $ s>0 $ we have
\begin{eqnarray*}
e^{-\frac{\tau}{2s}}\mbox{\boldmath $v$}_{\tau}
=(\mbox{\boldmath $e$}_{1}+{\rm i}\mbox{\boldmath $e$}_{2})\exp\left\{-\tau\left(\frac{\mbox{\boldmath $y$}-\mbox{\boldmath $x$}}{|\mbox{\boldmath $y$}-\mbox{\boldmath $x$}|^{2}}\cdot\mbox{\boldmath $e$}_2+\frac{1}{2s}\right)\right\}
\exp\left\{-{\rm i}\tau\frac{\mbox{\boldmath $y$}-\mbox{\boldmath $x$}}{|\mbox{\boldmath $y$}-\mbox{\boldmath $x$}|^{2}}\cdot\mbox{\boldmath $e$}_1\right\}.
\end{eqnarray*}
Since the real part of the power exponent takes the form
\begin{eqnarray}
\frac{\mbox{\boldmath $y$}-\mbox{\boldmath $x$}}{|\mbox{\boldmath $y$}-\mbox{\boldmath $x$}|^{2}}\cdot\mbox{\boldmath $e$}_{2}+\frac{1}{2s}
=\frac{2s(\mbox{\boldmath $y$}-\mbox{\boldmath $x$})\cdot\mbox{\boldmath $e$}_{2}+\vert \mbox{\boldmath $y$}-\mbox{\boldmath $x$}\vert^{2}}{2s\vert \mbox{\boldmath $y$}-\mbox{\boldmath $x$}\vert^{2}}
=\frac{\vert \mbox{\boldmath $y$}-(\mbox{\boldmath $x$}-s\mbox{\boldmath $e$}_{2})\vert^{2}-s^{2}}{2s\vert \mbox{\boldmath $y$}-\mbox{\boldmath $x$}\vert^{2}},
\label{eq:1-8}
\end{eqnarray}
one sees that $e^{-\frac{\tau}{2s}}\mbox{\boldmath $v$}_{\tau}$ has different asymptotic behaviors as $\tau\longrightarrow\infty$ whether $\mbox{\boldmath $y$}$ is inside or outside of the circle centred at $\mbox{\boldmath $x$}-s\mbox{\boldmath $e$}_2$ with radius $s$.
Namely, it holds the following facts;
\begin{itemize}
\item $\displaystyle \lim_{\tau\longrightarrow\infty}e^{-\frac{\tau}{2s}}\vert \mbox{\boldmath $v$}_{\tau}\vert=0 $ when $ \vert \mbox{\boldmath $y$}-(\mbox{\boldmath $x$}-s\mbox{\boldmath $e$}_{2})\vert>s $;
\item $\displaystyle \lim_{\tau\longrightarrow\infty}e^{-\frac{\tau}{2s}}\vert \mbox{\boldmath $v$}_{\tau}\vert=\infty $ when $\vert \mbox{\boldmath $y$}-(\mbox{\boldmath $x$}-s\mbox{\boldmath $e$}_{2})\vert<s $;
\item each components of $ \displaystyle e^{-\frac{\tau}{2s}}\mbox{\boldmath $v$}_{\tau} $ is highly oscillating as $ \tau\longrightarrow\infty $ when $ \vert \mbox{\boldmath $y$}-(\mbox{\boldmath $x$}-s\mbox{\boldmath $e$}_{2})\vert=s $.
\end{itemize}

Next, for fixed $\epsilon >0$ we put $ \mbox{\boldmath $x$}\in {\mathbb R}^2\setminus\overline{\Omega} $ on a line segment $\Gamma_{\epsilon}:= [0,a]\times\{b+\epsilon\}$.
Using the function $\mbox{\boldmath $v$}_{\tau}(\mbox{\boldmath $y$};\mbox{\boldmath $x$})$, we define a mathematical indicator and its derivative with respect to $\tau$.

\begin{definition}
Let $ \mbox{\boldmath $u$} $ be a weak solution of $ {\rm (*)} $. 
Given $ \mbox{\boldmath $x$}\in \Gamma_{\epsilon} $ and $ \tau >0 $ define
$$\left\{
\displaystyle
\begin{array}{l}
\displaystyle
I(\tau ; \mbox{\boldmath $x$})=
\int_{\partial\Omega}\mbox{\boldmath $g$}(\mbox{\boldmath $y$})\cdot \mbox{\boldmath $v$}_{\tau}(\mbox{\boldmath $y$}; \mbox{\boldmath $x$})-\mbox{\boldmath $u$}(\mbox{\boldmath $y$})\cdot \mbox{\boldmath $\sigma$}
(\mbox{\boldmath $v$}_{\tau}(\mbox{\boldmath $y$}; \mbox{\boldmath $x$}))\mbox{\boldmath $\nu$}\; {\rm d}s_{\mbox{\boldmath $y$}},
\\
\\
\displaystyle
I'(\tau ; \mbox{\boldmath $x$})=
\int_{\partial\Omega}\mbox{\boldmath $g$}(\mbox{\boldmath $y$})\cdot \mbox{\boldmath $v$}_{\tau}'(\mbox{\boldmath $y$}; \mbox{\boldmath $x$})-\mbox{\boldmath $u$}(\mbox{\boldmath $y$})\cdot \mbox{\boldmath $\sigma$}(\mbox{\boldmath $v$}_{\tau}'(\mbox{\boldmath $y$}; \mbox{\boldmath $x$}))\mbox{\boldmath $\nu$}\; {\rm d}s_{\mbox{\boldmath $y$}}.
\end{array}\right.
$$
Here and hereafter, $\mbox{\boldmath $v$}_{\tau}'=\mbox{\boldmath $v$}_{\tau}'(\mbox{\boldmath $y$}; \mbox{\boldmath $x$})$ denotes the derivative of $\mbox{\boldmath $v$}_{\tau}$ with respect to $\tau$.
\end{definition}

One of identification procedure for multiple cracks $\Sigma$ provided in \cite{IIS} is moving the virtual disc $B_{\tilde{s}}(\mbox{\boldmath $x$}-\tilde{s}\mbox{\boldmath $e$}_{2})$ with $\tilde{s}:= (b+\epsilon -c)/2$ while $ \mbox{\boldmath $x$} $ is varying on the line $\Gamma_{\epsilon}$, where we use the notation $ B_{s}(\mbox{\boldmath $x$}) := \{ \mbox{\boldmath $y$}\in {\mathbb R}^2 \; |\; |\mbox{\boldmath $y$} - \mbox{\boldmath $x$}| < s\} $ with $ s>0 $.
Now we extend Theorem 1 in \cite{IIS} to our Problem.

\begin{theorem}\label{th0}
Let $\mbox{\boldmath $g$}\in L^2(\partial \Omega)$ satisfy \eqref{eq:1-6} and the following condition {\rm ($\dag$)}:
\begin{itemize}
\item[{\rm ($\dag$)}] ${\rm supp}(\mbox{\boldmath $g$}) \subset ( \partial\Omega \cap \{|x_2-c|>\gamma \} ) \setminus \left(B_\gamma (\mbox{\boldmath $O$}) \cup B_\gamma (0,b) \cup B_\gamma (a,b) \cup B_\gamma (a,0) \right)$
for some $\gamma > 0$ and there exists $ \mbox{\boldmath $\rho$}_{0}\in {\cal R} $ such that either
$$\int_{\partial\Omega \cap \{x_2>c\} } \mbox{\boldmath $g$}\cdot \mbox{\boldmath $\rho$}_{0} \; {\rm d}s_{\mbox{\boldmath $x$}} \neq 0 \quad \text{or} \quad \int_{\partial\Omega \cap \{x_2<c\} } \mbox{\boldmath $g$}\cdot \mbox{\boldmath $\rho$}_{0} \; {\rm d}s_{\mbox{\boldmath $x$}} \neq 0,$$ 
\end{itemize}
Then, the following statements hold.
\begin{itemize}
\item[(S1)] if $\Gamma_{\epsilon}\ni\mbox{\boldmath $x$}=(x_1 , b+\epsilon) $ 
and $x_{1}\in \{c_{1},\cdots , c_{2m}\} $, then there exists an integer $N\geq 1$ and a complex number $\tilde{C}\neq 0$ such that 
\[
\lim_{\tau\longrightarrow\infty}\tau^{\frac{2N-1}{2}}
e^{-\frac{\tau}{2\tilde{s}}}I(\tau ; (x_{1}, b+\epsilon ))=\tilde{C}.
\]

\item[(S2)] if $ \mbox{\boldmath $x$}\in \Gamma_{\epsilon}\setminus (\{c_{1},\cdots , c_{2m}\}\times \{b+\epsilon\}) $, then $ e^{-\frac{\tau}{2\tilde{s}}}I(\tau ; (x_{1}, b+\epsilon )) $ is exponentially decaying as $\tau\longrightarrow\infty $.
\end{itemize}
\end{theorem}

Thus, in principle one can detect the location of all crack tips.
However, from numerical point of view, it may not be easy to catch the difference between algebraically and exponentially decaying because of the presence of measurements error and noise.
Therefore we consider another method used in \cite{HIIS, HI2020}, namely changing $s$.

Next, we define the function of $\mbox{\boldmath $x$}$ given by
\[
s_{\Sigma}(\mbox{\boldmath $x$}) := \sup\left\{s > 0\; |\; B_{s}( \mbox{\boldmath $x$}-s\mbox{\boldmath $e$}_{2})\subset {\mathbb R}^2\setminus\Sigma \right\}.
\]
Then one sees that the value $s_{\Sigma}(\mbox{\boldmath $x$})$ at $\mbox{\boldmath $x$}\in\Gamma_{\epsilon} $ implies the largest radius of the disc $ B_{s}( \mbox{\boldmath $x$}-s\mbox{\boldmath $e$}_{2}) $ whose exterior encloses $\Sigma$.
Now we describe the second result to extract the $s_{\Sigma}(\mbox{\boldmath $x$})$ at $\mbox{\boldmath $x$}\in\Gamma_{\epsilon}$ and a quantity specifying the position of the crack tips $c_1, \cdots , c_{2m}$.

\begin{theorem}\label{th1}
Let $\mbox{\boldmath $g$}\in L^2(\partial \Omega)$ satisfy \eqref{eq:1-6} and the condition {\rm ($\dag$)}.
Assume that $\mbox{\boldmath $x$}\in\Gamma_{\epsilon}$ satisfies the following condition {\rm ($\ddag$)}:
\begin{itemize}
\item[{\rm ($\ddag$)}] $\exists! j\in\{1,\cdots,2m\}$ s.t. $\displaystyle \overline{B_{s_{\Sigma}(\mbox{\boldmath $x$})}(\mbox{\boldmath $x$}-s_{\Sigma}(\mbox{\boldmath $x$})\mbox{\boldmath $e$}_2)}\cap\Sigma=\{(c_j,c)\}$.
\end{itemize}
Let the real number $\alpha \in] -\frac{\pi}{2},\, \frac{\pi}{2}]$ be the unique solution of the equation
\begin{eqnarray}
e^{{\rm i}\alpha} =(-1)^j\frac{c_{j}-x_{1}}{s_{\Sigma}(\mbox{\boldmath $x$})}+{\rm i}\frac{x_{2}-s_{\Sigma}(\mbox{\boldmath $x$})-c}{s_{\Sigma}(\mbox{\boldmath $x$})}.\label{eq:1-9}
\end{eqnarray}
Then, there exists a positive number $\tau_{0}$ such that for all $\tau\geq\tau_{0}$, $|I(\tau ; \mbox{\boldmath $x$})| > 0$, and we have the following formulae:
\begin{eqnarray}
\lim_{\tau\longrightarrow\infty}\frac{\log\vert I(\tau ;\mbox{\boldmath $x$})\vert}{\tau} &=& \frac{1}{2s_{\Sigma}(\mbox{\boldmath $x$})},\label{eq:1-10} \\
\lim_{\tau\longrightarrow\infty}\frac{I'(\tau ;\mbox{\boldmath $x$})}{I(\tau ;\mbox{\boldmath $x$})} &=& \frac{1}{2s_{\Sigma}(\mbox{\boldmath $x$})}+{\rm i}\frac{(-1)^{j+1}\cos{\alpha}}{2s_{\Sigma}(\mbox{\boldmath $x$})(1+\sin{\alpha})}.\label{eq:1-11}
\end{eqnarray}
\end{theorem}

\begin{figure}[!h]
\centering\includegraphics[width=8.5cm]{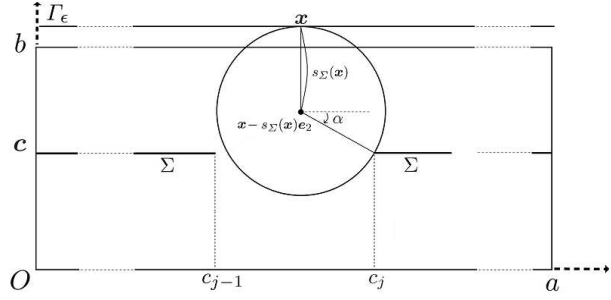}
\caption{an example of the angle $ \alpha $ in a case of even $j$}
\label{fig1-2}
\end{figure}

\begin{remark} $\;$
\begin{itemize}
\item For the meaning of $\alpha$, see Figure \ref{fig1-2}.
\item The formulae \eqref{eq:1-10} and \eqref{eq:1-11}
yield the value of $s_{\Sigma}(\mbox{\boldmath $x$})$.
Besides, taking the imaginary part of the both sides 
of formula \eqref{eq:1-11}, one gets
$$\displaystyle
(-1)^{j+1}\frac{\cos\alpha}{1+\sin\alpha}
=(-1)^{j+1}\tan\frac{1}{2}\left(\frac{\pi}{2}-\alpha\,\right).
$$
If the imaginary part of \eqref{eq:1-11} vanishes, then it means $\alpha =\frac{\pi}{2}$ which corresponds to the case of Theorem \ref{th0}, namely $s_{\Sigma}(\mbox{\boldmath $x$})=\tilde{s}$.
Otherwise, since $\vert\alpha\vert<\frac{\pi}{2}$, this enables us to decide whether $j$ is odd or even and then the value of $\alpha$ itself.
\item There are some examples of $ \mbox{\boldmath $g$}\in L^{2}(\partial\Omega) $ satisfying the conditions \eqref{eq:1-6} and {\rm ($\dag$)}.
In \cite{HI2020}, two concrete examples of such $ \mbox{\boldmath $g$} $ are provided, that is, for $\beta_{1}$, $\beta_{2} \neq 0$ and $\gamma' := \min\left\{ c-2\gamma , b-c-2\gamma \right\}$,
\[
\mbox{\boldmath $g$}_{1}=\left\{ \begin{array}{l}
\vspace{0.3cm}
\beta_{1}\mbox{\boldmath $e$}_{2} \quad {\rm on}\quad ] \gamma ,\, a-\gamma [\, \times \{b\},\\
\vspace{0.3cm}
-\beta_{1}\mbox{\boldmath $e$}_{2} \quad {\rm on}\quad ] \gamma ,\, a-\gamma [\, \times \{0\},\\
\mbox{\boldmath $0$}, \quad {\rm otherwise},
\end{array}\right.
\]
and
\[
\mbox{\boldmath $g$}_{2}=\left\{ \begin{array}{l}
\vspace{0.3cm}
\beta_{2} (a\mbox{\boldmath $e$}_{1}-(2\gamma +\gamma')\mbox{\boldmath $e$}_{2} ) \quad {\rm on}\quad \{a\}\times\, ] c-\gamma -\gamma' ,\, c-\gamma [,\\
\vspace{0.3cm}
-\beta_{2} (a\mbox{\boldmath $e$}_{1}-(2\gamma +\gamma')\mbox{\boldmath $e$}_{2} ) \quad {\rm on}\quad \{0\}\times \, ] c+\gamma ,\, c+\gamma +\gamma' [ ,\\
\mbox{\boldmath $0$},\quad {\rm otherwise}.
\end{array}\right.
\]
\item The condition {\rm ($\ddag$)} for $\mbox{\boldmath $x$}=(x_1,x_2)\in\Gamma_{\epsilon}$ implies that there exists a $j'\in\{1,2,\cdots ,m \}$ such that $x_{1}\in [c_{2j'-1}, c_{2j'}]\setminus \{(c_{2j'-1}+c_{2j'})/2 \} $ and $x_2 =b+\epsilon$.
In regard to a case that $\mbox{\boldmath $x$}\in\Gamma_{\epsilon}$ violates the condition {\rm ($\ddag$)} such as the projection of $\mbox{\boldmath $x$}=(x_1,x_2)\in\Gamma_{\epsilon}$ onto the line $x_{2}=c$ belongs to $\Sigma\setminus \{c_{1},\cdots , c_{2m}\}$, we should apply Theorem \ref{th0}.
\end{itemize}
\end{remark}

\section{Proof of Theorem \ref{th0} and \ref{th1}}

In order to prove Theorem \ref{th0} and \ref{th1} it is important to study the asymptotic behaviors of the indicator function and its derivative.
Firstly, using the Green formula, we obtain the following representation formulae of the indicator functions.

\begin{proposition}[Proposition 2 in \cite{II2007}]\label{prop1}
\begin{eqnarray}
I(\tau ;\mbox{\boldmath $x$}) &=& -\int_{\Sigma} \left( \mbox{\boldmath $u$}^{+}(\mbox{\boldmath $y$}) - \mbox{\boldmath $u$}^{-}(\mbox{\boldmath $y$})\right) \cdot \left( \mbox{\boldmath $\sigma$}(\mbox{\boldmath $v$}_{\tau}(\mbox{\boldmath $y$}; \mbox{\boldmath $x$}))\mbox{\boldmath $e$}_{2}\right)\; {\rm d}s_{\mbox{\boldmath $y$}}, \label{eq:2-3}\\
I'(\tau ;\mbox{\boldmath $x$}) &=& -\int_{\Sigma} \left( \mbox{\boldmath $u$}^{+}(\mbox{\boldmath $y$}) - \mbox{\boldmath $u$}^{-}(\mbox{\boldmath $y$})\right) \cdot \left( \mbox{\boldmath $\sigma$}(\mbox{\boldmath $v$}_{\tau}'(\mbox{\boldmath $y$}; \mbox{\boldmath $x$}))\mbox{\boldmath $e$}_{2}\right)\; {\rm d}s_{\mbox{\boldmath $y$}}. \label{eq:2-31}
\end{eqnarray}
\end{proposition}
\noindent Here $\mbox{\boldmath $u$}^{\pm}$ denotes the trace of $\mbox{\boldmath $u$}|_{\Omega^{\pm}}$ onto $]0,a[ \times \{c\}$, respectively.
By use of \eqref{eq:1-3} and \eqref{eq:1-1} one has
\begin{eqnarray}
\mbox{\boldmath $\sigma$}(\mbox{\boldmath $v$}_{\tau}(\mbox{\boldmath $y$}; \mbox{\boldmath $x$}))\mbox{\boldmath $e$}_{2} &=& 2\mu\tau \left( \frac{\mbox{\boldmath $y$}-\mbox{\boldmath $x$}}{|\mbox{\boldmath $y$}-\mbox{\boldmath $x$}|^{2}}\cdot ( \mbox{\boldmath $e$}_{2} + {\rm i} \mbox{\boldmath $e$}_{1})\right)^{2}\mbox{\boldmath $v$}_{\tau}, \label{eq:2-7} \\
\mbox{\boldmath $\sigma$}(\mbox{\boldmath $v$}_{\tau}'(\mbox{\boldmath $y$}; \mbox{\boldmath $x$}))\mbox{\boldmath $e$}_{2} &=& 2\mu \left( 1-\tau \frac{\mbox{\boldmath $y$}-\mbox{\boldmath $x$}}{|\mbox{\boldmath $y$}-\mbox{\boldmath $x$}|^{2}}\cdot ( \mbox{\boldmath $e$}_{2} + {\rm i} \mbox{\boldmath $e$}_{1})\right)\left( \frac{\mbox{\boldmath $y$}-\mbox{\boldmath $x$}}{|\mbox{\boldmath $y$}-\mbox{\boldmath $x$}|^{2}}\cdot ( \mbox{\boldmath $e$}_{2} + {\rm i} \mbox{\boldmath $e$}_{1})\right)^{2}\mbox{\boldmath $v$}_{\tau}. \label{eq:2-8}
\end{eqnarray}

\subsection{Proof of Theorem \ref{th0}}

(S1) can be deduced from a combination of Proposition \ref{prop2} and Lemma \ref{lem1} described later, the same as (S1) of \cite[Theorem 1]{IIS}.
Thus, it suffices to prove only (S2).
The proof proceeds along the same line as in \cite[Section 2.2]{IIS}.
Indeed, we firstly consider the case when $\Gamma_{\epsilon}\ni\mbox{\boldmath $x$}=(x_1 , b+\epsilon) $ and $x_{1}\in\bigcup_{j=1}^{m}] c_{2j-1}, c_{2j}[$.
Since $ \vert \mbox{\boldmath $y$}-(\mbox{\boldmath $x$}-\tilde{s}\mbox{\boldmath $e$}_{2})\vert^2 -\tilde{s}^{2}\geq \min_{i=1,\cdots, 2m}|x_{1}-c_{i}|^2 >0 $ for all $\mbox{\boldmath $y$}\in\Sigma $, it follows from \eqref{eq:1-8}, \eqref{eq:2-3} and \eqref{eq:2-7} that there exists a $\alpha_{0}>0$ such that
\[
\left| e^{-\frac{\tau}{2\tilde{s}}}I(\tau ; (x_{1}, b+\epsilon )) \right| \leq \left\| \mbox{\boldmath $u$}^{+} - \mbox{\boldmath $u$}^{-}\right\|_{L^{2}(\Sigma)} \cdot \left\| \mbox{\boldmath $\sigma$}(\mbox{\boldmath $v$}_{\tau})\mbox{\boldmath $e$}_{2}\right\|_{L^{2}(\Sigma)}=O\left(\tau e^{-\alpha_{0}\tau}\right)
\]
as $\tau\longrightarrow\infty$.

Next, we consider the case when $\Gamma_{\epsilon}\ni\mbox{\boldmath $x$}=(x_1 , b+\epsilon) $ and $x_{1}\in\bigcup_{j=0}^{m}] c_{2j}, c_{2j+1}[$.
By use of Airy's stress function $U$, a problem \eqref{eq:1-4} in a neighborhood of $\mbox{\boldmath $x$}-2\tilde{s}\mbox{\boldmath $e$}_2$ and $x_2 >c$ with a traction free condition on $x_2=c$ can be reduced to a problem for biharmonic equation $\triangle^{2}U=0 $ with Dirichlet boundary condition $U=\partial U / \partial\mbox{\boldmath $\nu$}=0 $ (cf. \cite{II2009}).
According to an extension formula for $U$ across a straight line, refer to \cite{D1956, F1994}, one sees that in a neighborhood of $\mbox{\boldmath $x$}-2\tilde{s}\mbox{\boldmath $e$}_2$, $\mbox{\boldmath $u$}^{+}$ has a continuation $ \widetilde{\mbox{\boldmath $u$}^{+}} $ from $x_{2}>c$ into $x_{2}<c$.
Then we have that for a sufficiently small $\delta_{0}$ such as $\{c_{0},\cdots , c_{2m+1}\}\times \{ c\} \subset \Sigma\setminus\overline{B_{\tilde{s}+\delta_{0}}(\mbox{\boldmath $x$}-\tilde{s}\mbox{\boldmath $e$}_2)}$
\begin{eqnarray*}
\lefteqn{\hspace{1cm}
e^{-\frac{\tau}{2\tilde{s}}}\int_{\Sigma\cap B_{\tilde{s}+\delta_{0}}(\mbox{\boldmath $x$}-\tilde{s}\mbox{\boldmath $e$}_2)}\left( \mbox{\boldmath $u$}^{+}(\mbox{\boldmath $y$}) - \mbox{\boldmath $u$}^{-}(\mbox{\boldmath $y$})\right) \cdot \left( \mbox{\boldmath $\sigma$}(\mbox{\boldmath $v$}_{\tau}(\mbox{\boldmath $y$}; \mbox{\boldmath $x$}))\mbox{\boldmath $e$}_{2}\right)\; {\rm d}s_{\mbox{\boldmath $y$}}}\\
&&=e^{-\frac{\tau}{2\tilde{s}}}\int_{\partial B_{\tilde{s}+\delta_{0}}(\mbox{\boldmath $x$}-\tilde{s}\mbox{\boldmath $e$}_2)\cap\{x_{2}<c\}}
\left( \widetilde{\mbox{\boldmath $u$}^{+}}-\mbox{\boldmath $u$}^{-} \right)
\cdot \mbox{\boldmath $\sigma$}
(\mbox{\boldmath $v$}_{\tau})\mbox{\boldmath $\nu$}
-\mbox{\boldmath $v$}_{\tau}\cdot\mbox{\boldmath $\sigma$}\left( \tilde{\mbox{\boldmath $u$}^{+}}-\mbox{\boldmath $u$}^{-} \right)\mbox{\boldmath $\nu$}\; {\rm d}s_{\mbox{\boldmath $y$}},
\end{eqnarray*}
where $\mbox{\boldmath $\nu$}$ is the inward normal to $B_{\tilde{s}+\delta_{0}}(\mbox{\boldmath $x$}-\tilde{s}\mbox{\boldmath $e$}_2)$.
There exists a $\alpha_{1}>0$ such that the right-hand side has the bound $O\left(\tau e^{-\alpha_{1}\tau}\right)$ as $\tau\longrightarrow\infty$.
Since $\vert \mbox{\boldmath $y$}- (\mbox{\boldmath $x$}-\tilde{s}\mbox{\boldmath $e$}_2) \vert > \tilde{s}+\delta_{0} $ for $ \mbox{\boldmath $y$}\in \Sigma\setminus\overline{B_{\tilde{s}+\delta_{0}}(\mbox{\boldmath $x$}-\tilde{s}\mbox{\boldmath $e$}_2)} $, we conclude that
\[
 e^{-\frac{\tau}{2\tilde{s}}}I(\tau ; (x_{1}, b+\epsilon )) =O\left(\tau e^{-\alpha_{1}\tau}\right)
\]
as $\tau\longrightarrow\infty$.

\subsection{Proof of Theorem \ref{th1}}

For $\mbox{\boldmath $x$}\in\Gamma_{\epsilon}$ satisfying {\rm ($\ddag$)}, $(c_j , c)$ is uniquely determined and then we choose a polar coordinates system with respect to the center $(c_j , c)$ depending on whether $j$ is even or odd.
Let $\eta_{0}$ be a positive number such that
\[
\eta_{0} < \min_{i=1,\cdots , 2m+1}\left\{ c_{i}-c_{i-1}\right\} \quad {\rm and}\quad \eta_{0} < \min\{ b-c , c \}.
\]

\begin{itemize}
\item If $j$ is odd, then we set $ \mbox{\boldmath $x$}=(c_{j}+r\cos{\theta} , c+r\sin{\theta}) $ for $r\in \, ] 0, \eta_{0} [ $ and $\theta \in \, ]-\pi , \pi [ $, and define
\begin{eqnarray*}
\left\{ \begin{array}{ll}
\vspace{0.3cm}
\mbox{\boldmath $u$}^{+}(\mbox{\boldmath $x$})=\mbox{\boldmath $u$}(r,\theta ):=\mbox{\boldmath $u$} (c_{j}+r\cos{\theta} , c+r\sin{\theta}), & r\in \, ] 0, \eta_{0} [ , \quad \theta \in \, ]0, \pi [,\\
\mbox{\boldmath $u$}^{-}(\mbox{\boldmath $x$})=\mbox{\boldmath $u$}(r,\theta ):=\mbox{\boldmath $u$} (c_{j}+r\cos{\theta} , c+r\sin{\theta}), & r\in \, ] 0, \eta_{0} [ , \quad \theta \in \, ]-\pi, 0 [.
\end{array}\right.
\end{eqnarray*}
\item If $j$ is even, then we set $ \mbox{\boldmath $x$}=(c_{j}-r\cos{\theta} , c-r\sin{\theta}) $ for $r\in \, ] 0, \eta_{0} [ $ and $\theta \in \, ]-\pi , \pi [ $, and define
\begin{eqnarray*}
\left\{ \begin{array}{ll}
\vspace{0.3cm}
\mbox{\boldmath $u$}^{+}(\mbox{\boldmath $x$})=\mbox{\boldmath $u$}(r,\theta ):=\mbox{\boldmath $u$} (c_{j}-r\cos{\theta} , c-r\sin{\theta}), & r\in \, ] 0, \eta_{0} [ , \quad \theta \in \, ] -\pi , 0 [,\\
\mbox{\boldmath $u$}^{-}(\mbox{\boldmath $x$})=\mbox{\boldmath $u$}(r,\theta ):=\mbox{\boldmath $u$} (c_{j}-r\cos{\theta} , c-r\sin{\theta}), & r\in \, ] 0, \eta_{0} [ , \quad \theta \in \, ] 0, \pi [.
\end{array}\right.
\end{eqnarray*}
\end{itemize}

Next, we recall a convergent series expansion of a weak solution of $(*)$ around a tip of crack $(c_j,c)$.

\begin{proposition}[Proposition 1 in \cite{II2007, II2008}]\label{prop0}
Fix $\eta\in \, ]0, \frac{\eta_{0}}{2}[$.
There exist real numbers $ A_{k}^{(j)} $, $ B_{k}^{(j)} $ $ (k=0, 1, 2, \cdots ) $ and $\mbox{\boldmath $\rho$}_{j}\in {\cal R} $ such that
\begin{eqnarray}
\mbox{\boldmath $u$}(r,\theta ) -\mbox{\boldmath $\rho$}_{j} =\sum_{k=0}^{\infty}\frac{A_{k}^{(j)}}{2\mu}r^{\frac{k}{2}}\mbox{\boldmath $\varphi$}_{k}(\theta )-\sum_{k=0}^{\infty}\frac{B_{k}^{(j)}}{2\mu}r^{\frac{k}{2}}\mbox{\boldmath $\psi$}_{k}(\theta )\quad {\rm in}\quad B_{2\eta}((c_j,c))\setminus\Sigma , \label{eq:2-4}
\end{eqnarray}
where
\begin{eqnarray*}
\mbox{\boldmath $\varphi$}_{k}(\theta ) &=& \left( \begin{array}{c}
\vspace{0.3cm}
\kappa\cos{\frac{k}{2}\theta }-\frac{k}{2}\cos{\left(\frac{k}{2}-2\right)\theta }+\left\{\frac{k}{2}+(-1)^{k}\right\}\cos{\frac{k}{2}\theta } \\
\kappa\sin{\frac{k}{2}\theta }+\frac{k}{2}\sin{\left(\frac{k}{2}-2\right)\theta }-\left\{\frac{k}{2}+(-1)^{k}\right\}\sin{\frac{k}{2}\theta }
\end{array}\right), \\
\\
\mbox{\boldmath $\psi$}_{k}(\theta ) &=& \left( \begin{array}{c}
\vspace{0.3cm}
\kappa\sin{\frac{k}{2}\theta }-\frac{k}{2}\sin{\left(\frac{k}{2}-2\right)\theta }+\left\{\frac{k}{2}-(-1)^{k}\right\}\sin{\frac{k}{2}\theta } \\
-\kappa\cos{\frac{k}{2}\theta }-\frac{k}{2}\cos{\left(\frac{k}{2}-2\right)\theta }+\left\{\frac{k}{2}-(-1)^{k}\right\}\cos{\frac{k}{2}\theta }
\end{array}\right)
\end{eqnarray*}
with $ \kappa=\frac{\lambda +3\mu }{\lambda +\mu} $.
The series is convergent, absolutely in $ H^{1}(B_{\eta}((c_j,c))\cap \Omega^{+}) $ and $ H^{1}(B_{\eta}((c_j,c))\cap \Omega^{-}) $, and uniformly on compact sets in $ B_{2\eta}((c_j,c)) $.
Moreover, for each $ n=1,2, \cdots $, the following estimate is valid uniformly for $r\in \, ]0, \eta [ $
\begin{eqnarray*}
\lefteqn{\hspace{1cm}
\left| \mbox{\boldmath $u$}(r, \pi) -\mbox{\boldmath $\rho$}_{j} -\left( \sum_{k=1}^{n}\frac{A_{k}^{(j)}}{2\mu}r^{\frac{k}{2}}\mbox{\boldmath $\varphi$}_{k}(\pi )-\sum_{k=1}^{n}\frac{B_{k}^{(j)}}{2\mu}r^{\frac{k}{2}}\mbox{\boldmath $\psi$}_{k}(\pi )\right) \right|} \label{eq:2-5} \\
&&+\left| \mbox{\boldmath $u$}(r, -\pi) -\mbox{\boldmath $\rho$}_{j} -\left( \sum_{k=1}^{n}\frac{A_{k}^{(j)}}{2\mu}r^{\frac{k}{2}}\mbox{\boldmath $\varphi$}_{k}(-\pi )-\sum_{k=1}^{n}\frac{B_{k}^{(j)}}{2\mu}r^{\frac{k}{2}}\mbox{\boldmath $\psi$}_{k}(-\pi )\right) \right|
\leq K_{n} r^{\frac{n+1}{2}}, 
\end{eqnarray*}
where $K_n$ is a positive constant depending on $n$.
\end{proposition}

Then, by virtue of Proposition \ref{prop0} one sees that for each $n =1,2,\cdots$ and $ \mbox{\boldmath $y$}\in B_{\eta}((c_j,c))\cap \Sigma $
\begin{eqnarray}
\mbox{\boldmath $u$}^{+}(\mbox{\boldmath $y$}) - \mbox{\boldmath $u$}^{-}(\mbox{\boldmath $y$}) &=& (-1)^{j}\left( \mbox{\boldmath $u$}(r, -\pi ) - \mbox{\boldmath $u$}(r, \pi ) \right) \nonumber \\
&=&\frac{\kappa +1}{\mu} \sum_{k=1}^{n}(-1)^{k+j}r^{\frac{2k-1}{2}}\left( \begin{array}{c}
\vspace{0.3cm}
-B_{2k-1}^{(j)} \\
A_{2k-1}^{(j)} \end{array}\right) 
+O\left(r^{\frac{2n+1}{2}}\right) \mbox{\boldmath $k$}
\label{eq:2-6}
\end{eqnarray}
with a constant vector $\mbox{\boldmath $k$}$.

From \eqref{eq:2-7}, \eqref{eq:2-8} and \eqref{eq:2-6}, the following asymptotic expansions for $I(\tau ;\mbox{\boldmath $x$})$ and $I'(\tau ;\mbox{\boldmath $x$})$ are obtained.
\begin{proposition}\label{prop2}
As $\tau\longrightarrow\infty $, for each fixed $\mbox{\boldmath $x$}\in\Gamma_{\epsilon}$ satisfying {\rm ($\ddag$)}
\begin{eqnarray}
e^{-\frac{\tau}{2s_{\Sigma}(\mbox{\boldmath $x$})}}I(\tau ;\mbox{\boldmath $x$}) &=& \sum_{k=1}^{n}C_{k}e^{\frac{(-1)^{j+1}{\rm i}\tau\cos{\alpha}}{2s_{\Sigma}(\mbox{\boldmath $x$})(1+\sin{\alpha})}} \tau^{-\frac{2k-1}{2}}+O\left(\tau^{-\frac{n+1}{2}}\tau^{\frac{3}{4}}\right), \label{eq:2-9} \\
e^{-\frac{\tau}{2s_{\Sigma}(\mbox{\boldmath $x$})}}I'(\tau ;\mbox{\boldmath $x$}) &=& \sum_{k=1}^{n}C_{k}\left( 1+\tau\left( \frac{1}{2s_{\Sigma}(\mbox{\boldmath $x$})}+\frac{(-1)^{j+1}{\rm i}\cos{\alpha}}{2s_{\Sigma}(\mbox{\boldmath $x$})(1+\sin{\alpha})} \right) \right) e^{\frac{(-1)^{j+1}{\rm i}\tau\cos{\alpha}}{2s_{\Sigma}(\mbox{\boldmath $x$})(1+\sin{\alpha})}} \tau^{-\frac{2k+1}{2}}\nonumber \\
&&+O\left(\tau^{-\frac{n+1}{2}}\tau^{\frac{3}{4}}\right), \label{eq:2-10} 
\end{eqnarray}
where
\begin{eqnarray*}
C_{k}&:=&{\rm i}(\kappa +1)(-1)^{k}\left( s_{\Sigma}(\mbox{\boldmath $x$})\right)^{2k-1}2^{\frac{2k+1}{2}}(1+\sin{\alpha})^{\frac{2k-1}{2}}e^{(-1)^{j+1}{\rm i}\left(\frac{2k-1}{2}\right)\alpha}\Gamma\left(\frac{2k+1}{2} \right)\\
&& \times
\left( \begin{array}{c}
\vspace{0.3cm}
B_{2k-1}^{(j)} \\
-A_{2k-1}^{(j)} \end{array}\right)\cdot (\mbox{\boldmath $e$}_{1}+{\rm i}\mbox{\boldmath $e$}_{2})
\end{eqnarray*}
\end{proposition}

The proof of Proposition \ref{prop2} is given in Section 5, however this is not sufficient to prove Theorem \ref{th0} and \ref{th1} because it only shows that $e^{-\tau/(2s_{\Sigma}(\mbox{\boldmath $x$}))}I(\tau ;\mbox{\boldmath $x$})$ and $e^{-\tau/(2s_{\Sigma}(\mbox{\boldmath $x$}))}I'(\tau ;\mbox{\boldmath $x$})$ are {\it at most} algebraically decaying as $\tau\longrightarrow\infty$.
In order to complete the proof of Theorem \ref{th0} and \ref{th1}, we need to investigate non-vanishing of a coefficient in the expansion \eqref{eq:2-9} and \eqref{eq:2-10}, noting that $ C_{n}=0 $ if and only if $ A_{2n-1}^{(j)}=B_{2n-1}^{(j)}=0$.

\begin{lemma}\label{lem1}
Let $ \mbox{\boldmath $g$} $ satisfy the conditions \eqref{eq:1-6} and {\rm ($\dag$)}.
Then, there exists an integer $n\geq 1$ such that $\left(A_{2n-1}^{(j)}\right)^{2}+\left(B_{2n-1}^{(j)}\right)^{2}\neq 0$.
\end{lemma}

\noindent The proof of Lemma \ref{lem1} is given in Section 5, which is not given a detailed description in \cite{HI2020}.

In the consequence we can take $N:= \min\left\{n\geq 1\; \left|\; \left(A_{2n-1}^{(j)}\right)^{2}+\left(B_{2n-1}^{(j)}\right)^{2}\neq 0\right. \right\}$.
Then, substituting $n=2N$ in Proposition \ref{prop2}, we obtain that for each $ \mbox{\boldmath $x$}\in\Gamma_{\epsilon} $ satisfying the condition {\rm ($\ddag$)}
\begin{eqnarray*}
\lim_{\tau\longrightarrow\infty}\tau^{\frac{2N-1}{2}}e^{\frac{(-1)^{j}{\rm i}\tau\cos{\alpha}}{2s_{\Sigma}(\mbox{\boldmath $x$})(1+\sin{\alpha})}}e^{-\frac{\tau}{2s_{\Sigma}(\mbox{\boldmath $x$})}}I(\tau ;\mbox{\boldmath $x$})&=&C_{N}\neq 0, \\
\lim_{\tau\longrightarrow\infty}\tau^{\frac{2N-1}{2}}e^{\frac{(-1)^{j}{\rm i}\tau\cos{\alpha}}{2s_{\Sigma}(\mbox{\boldmath $x$})(1+\sin{\alpha})}}e^{-\frac{\tau}{2s_{\Sigma}(\mbox{\boldmath $x$})}}I'(\tau ;\mbox{\boldmath $x$})&=&C_{N}\left(\frac{1}{2s_{\Sigma}(\mbox{\boldmath $x$})}+{\rm i}\frac{(-1)^{j+1}\cos{\alpha}}{2s_{\Sigma}(\mbox{\boldmath $x$})(1+\sin{\alpha})}\right),
\end{eqnarray*}
which immediately yield (S1) in Theorem \ref{th0} and Theorem \ref{th1}.

\section{Proof of Proposition \ref{prop2} and Lemma \ref{lem1}}
\subsection{Proof of Proposition \ref{prop2}}
The proof of \eqref{eq:2-9} type formula of Proposition \ref{prop2}
in the case of the Laplace equation is given in \cite{HIIS}.
However \eqref{eq:2-10} type formula is not considered therein.
Some calculation such as estimates of oscillating integrals in \cite{HIIS} works also for the present case, therefore we just referred them in this section.

We fix $\mbox{\boldmath $x$}\in\Gamma_{\epsilon}$ satisfying {\rm ($\ddag$)} and in what follows use the notation for simplicity $s_{0}:=s_{\Sigma}(\mbox{\boldmath $x$})$.
Choose $\delta >0$ in such a way that
\begin{eqnarray*}
\overline{B_{s_0+\delta}(\mbox{\boldmath $x$}-s_0\mbox{\boldmath $e$}_2)}\cap\Sigma \subset [c_{j-1},c_j]\times\{c\} & {\rm if} & j\;\; {\rm is}\;\;{\rm odd},\\
\overline{B_{s_0+\delta}(\mbox{\boldmath $x$}-s_0\mbox{\boldmath $e$}_2)}\cap\Sigma \subset [c_{j},c_{j+1}]\times\{c\} & {\rm if} & j\;\; {\rm is}\;\; {\rm even},
\end{eqnarray*}
and
\[
\eta_{\delta}:=\sqrt{(s_0+\delta)^{2}-s_0^2}<\eta 
\]
with $\eta$ fixed in Proposition \ref{prop0}.
Set 
\[
\eta_{\delta}':=\sqrt{(s_{0}+\delta)^2-(x_{2}-s_{0}-c)^2}-\vert x_{1}-c_j\vert ,
\]
and then it is easy to see $\eta_{\delta}>\eta'_{\delta}>0$.
Now we divide $\Sigma$ into two disjoint parts:
\[
\Sigma = \left( \Sigma\setminus B_{s_0+\delta}(\mbox{\boldmath $x$}-s_0\mbox{\boldmath $e$}_2)\right) \cup \left( \Sigma\cap B_{s_0+\delta}(\mbox{\boldmath $x$}-s_0\mbox{\boldmath $e$}_2)\right).
\]

For $\mbox{\boldmath $y$}\in \Sigma\setminus B_{s_0+\delta}(\mbox{\boldmath $x$}-s_0\mbox{\boldmath $e$}_2)$, it holds $| \mbox{\boldmath $y$}-(\mbox{\boldmath $x$}-s_{0}\mbox{\boldmath $e$}_{2}) |\geq s_{0}+\delta $.
Consequently, it follows from \eqref{eq:1-8}, \eqref{eq:2-7} and \eqref{eq:2-8} that for all $ \mbox{\boldmath $y$}\in \Sigma\setminus B_{s_0+\delta}(\mbox{\boldmath $x$}-s_0\mbox{\boldmath $e$}_2) $
\[
e^{-\frac{\tau}{2s_{0}}}\left( \left| \mbox{\boldmath $\sigma$}(\mbox{\boldmath $v$}_{\tau}(\mbox{\boldmath $y$}; \mbox{\boldmath $x$}))\mbox{\boldmath $e$}_{2} \right|+ \left| \mbox{\boldmath $\sigma$}(\mbox{\boldmath $v$}_{\tau}'(\mbox{\boldmath $y$}; \mbox{\boldmath $x$}))\mbox{\boldmath $e$}_{2} \right|\right) \leq \alpha_{2}(1+ \alpha_{3}\tau ) e^{-\alpha_{4}\eta_{\delta}^{2}\tau}
\]
with positive constants $\alpha_{2}$, $\alpha_{3}$ and $\alpha_{4}$.
Combining with \eqref{eq:2-6}, one can see that
\begin{eqnarray}
e^{-\frac{\tau}{2s_{0}}}\int_{\Sigma\setminus B_{s_0+\delta}(\mbox{\boldmath $x$}-s_0\mbox{\boldmath $e$}_2)}\left( \mbox{\boldmath $u$}^{+}(\mbox{\boldmath $y$}) - \mbox{\boldmath $u$}^{-}(\mbox{\boldmath $y$})\right) \cdot \left( \mbox{\boldmath $\sigma$}(\mbox{\boldmath $v$}_{\tau}(\mbox{\boldmath $y$}; \mbox{\boldmath $x$}))\mbox{\boldmath $e$}_{2}\right)\, {\rm d}s_{\mbox{\boldmath $y$}} &=& O(\tau e^{-\alpha_{4}\eta_{\delta}^{2}\tau}),
\label{eq:3-5} \\
e^{-\frac{\tau}{2s_{0}}}\int_{\Sigma\setminus B_{s_0+\delta}(\mbox{\boldmath $x$}-s_0\mbox{\boldmath $e$}_2)}\left( \mbox{\boldmath $u$}^{+}(\mbox{\boldmath $y$}) - \mbox{\boldmath $u$}^{-}(\mbox{\boldmath $y$})\right) \cdot \left( \mbox{\boldmath $\sigma$}(\mbox{\boldmath $v$}_{\tau}'(\mbox{\boldmath $y$}; \mbox{\boldmath $x$}))\mbox{\boldmath $e$}_{2}\right)\, {\rm d}s_{\mbox{\boldmath $y$}} &=& O(\tau e^{-\alpha_{4}\eta_{\delta}^{2}\tau}),
\label{eq:3-51}
\end{eqnarray}
as $\tau \longrightarrow \infty$.\\

Next, for $\mbox{\boldmath $y$}\in \Sigma\cap B_{s_0+\delta}(\mbox{\boldmath $x$}-s_0\mbox{\boldmath $e$}_2)$, noting \eqref{eq:2-6}, \eqref{eq:2-7} and Lemma 3.3 in \cite{HIIS}, we have
\begin{eqnarray}
\lefteqn{\hspace{1cm}
-e^{-\frac{\tau}{2s_{0}}}\int_{\Sigma\cap B_{s_0+\delta}(\mbox{\boldmath $x$}-s_0\mbox{\boldmath $e$}_2)}\left( \mbox{\boldmath $u$}^{+}(\mbox{\boldmath $y$}) - \mbox{\boldmath $u$}^{-}(\mbox{\boldmath $y$})\right) \cdot \left( \mbox{\boldmath $\sigma$}(\mbox{\boldmath $v$}_{\tau}(\mbox{\boldmath $y$}; \mbox{\boldmath $x$}))\mbox{\boldmath $e$}_{2}\right)\, {\rm d}s_{\mbox{\boldmath $y$}}} \nonumber \\
&& =-e^{-\frac{\tau}{2s_{0}}} \int_{\Sigma\cap B_{s_0+\delta}(\mbox{\boldmath $x$}-s_0\mbox{\boldmath $e$}_2 ) }
\left\{ \frac{\kappa +1}{\mu} \sum_{k=1}^{n}(-1)^{k+j}r^{\frac{2k-1}{2}}\left( \begin{array}{c}
\vspace{0.3cm}
-B_{2k-1}^{(j)} \\
A_{2k-1}^{(j)} \end{array}\right) +O\left(r^{\frac{2n+1}{2}}\right) \mbox{\boldmath $k$} \right\} \nonumber \\
&& \quad \cdot \left\{ 2\mu\tau \left( \frac{\mbox{\boldmath $y$}-\mbox{\boldmath $x$}}{|\mbox{\boldmath $y$}-\mbox{\boldmath $x$}|^{2}}\cdot ( \mbox{\boldmath $e$}_{2} + {\rm i} \mbox{\boldmath $e$}_{1})\right)^{2}\exp{\left\{ -\tau \frac{\mbox{\boldmath $y$}-\mbox{\boldmath $x$}}{|\mbox{\boldmath $y$}-\mbox{\boldmath $x$}|^{2}}\cdot (\mbox{\boldmath $e$}_{2}+{\rm i}\mbox{\boldmath $e$}_{1}) \right\}}(\mbox{\boldmath $e$}_{1}+{\rm i}\mbox{\boldmath $e$}_{2}) \right\}\; {\rm d}s_{\mbox{\boldmath $y$}} \nonumber \\
&& = 2(\kappa +1)\tau e^{-\frac{\tau}{2s_{0}}}\sum_{k=1}^{n}(-1)^{k+j}\left( \begin{array}{c}
\vspace{0.3cm}
B_{2k-1}^{(j)} \\
-A_{2k-1}^{(j)} \end{array}\right)\cdot (\mbox{\boldmath $e$}_{1}+{\rm i}\mbox{\boldmath $e$}_{2})I_{k}(\tau)+ O\left(\tau^{-\frac{n+1}{2}}\tau^{\frac{3}{4}}\right).
\label{eq:3-6}
\end{eqnarray}

\noindent Here $ I_{k}(\tau) $ depending on $j$ is defined as follows:
\begin{eqnarray*}
I_{k}(\tau ) := \left\{ \begin{array}{ccc}
\vspace{0.3cm}
\displaystyle \int_{0}^{\eta_{\delta}'}\frac{r^{\frac{2k-1}{2}}}{(r-s_{0}\overline{z_{\alpha}})^{2}}\exp\left( \frac{{\rm i}\tau}{r-s_{0}\overline{z_{\alpha}}} \right)\; {\rm d}r & {\rm if} & j\;\; {\rm is}\;\;{\rm odd},\\
\displaystyle \int_{0}^{\eta_{\delta}'}\frac{r^{\frac{2k-1}{2}}}{(r-s_{0}z_{\alpha})^{2}}\exp\left( \frac{-{\rm i}\tau}{r-s_{0}z_{\alpha}} \right)\; {\rm d}r & {\rm if} & j\;\; {\rm is}\;\;{\rm even},
\end{array} \right.
\end{eqnarray*}
where
\begin{eqnarray*}
z_{\alpha} := -\left( e^{-\frac{\pi}{2}{\rm i}}+{\rm i}e^{-(\frac{\pi}{2}+\alpha ){\rm i}} \right)
=-\cos{\alpha} +{\rm i} (1+\sin{\alpha}),
\end{eqnarray*}
$\alpha \in] -\frac{\pi}{2},\, \frac{\pi}{2}]$ is uniquely determined as the solution of \eqref{eq:1-9} and $\overline{z_{\alpha}}$ denotes the complex conjugate of $z_{\alpha}$.
In this notation for $ \mbox{\boldmath $y$}=(c_{j}+(-1)^{j}r , c) $ one sees
\[
\frac{\mbox{\boldmath $y$}-\mbox{\boldmath $x$}}{|\mbox{\boldmath $y$}-\mbox{\boldmath $x$}|^{2}}\cdot ( \mbox{\boldmath $e$}_{2} + {\rm i} \mbox{\boldmath $e$}_{1})=\frac{-{\rm i}}{(-1)^{j+1}r+(x_{1}-c_{j})+{\rm i}(x_{2}-c)}=\left\{ \begin{array}{ccc}
\vspace{0.3cm}
\displaystyle \frac{-{\rm i}}{r-s_{0}\overline{z_{\alpha}}} & {\rm if} & j\;\; {\rm is}\;\;{\rm odd},\\
\displaystyle \frac{{\rm i}}{r-s_{0}z_{\alpha}} & {\rm if} & j\;\; {\rm is}\;\;{\rm even}.
\end{array} \right.
\]
Similarly, for $\mbox{\boldmath $y$}\in \Sigma\cap B_{s_0+\delta}(\mbox{\boldmath $x$}-s_0\mbox{\boldmath $e$}_2)$, it follows from \eqref{eq:2-8} that
\begin{eqnarray}
\lefteqn{\hspace{1cm}
-e^{-\frac{\tau}{2s_{0}}}\int_{\Sigma\cap B_{s_0+\delta}(\mbox{\boldmath $x$}-s_0\mbox{\boldmath $e$}_2)}\left( \mbox{\boldmath $u$}^{+}(\mbox{\boldmath $y$}) - \mbox{\boldmath $u$}^{-}(\mbox{\boldmath $y$})\right) \cdot \left( \mbox{\boldmath $\sigma$}(\mbox{\boldmath $v$}'_{\tau}(\mbox{\boldmath $y$}; \mbox{\boldmath $x$}))\mbox{\boldmath $e$}_{2}\right)\, {\rm d}s_{\mbox{\boldmath $y$}}} \nonumber \\
&& = 2(\kappa +1) e^{-\frac{\tau}{2s_{0}}}\sum_{k=1}^{n}(-1)^{k+j}\left( \begin{array}{c}
\vspace{0.3cm}
B_{2k-1}^{(j)} \\
-A_{2k-1}^{(j)} \end{array}\right)\cdot (\mbox{\boldmath $e$}_{1}+{\rm i}\mbox{\boldmath $e$}_{2})I'_{k}(\tau)+ O\left(\tau^{-\frac{n+1}{2}}\tau^{\frac{3}{4}}\right),\nonumber \\
&&\qquad 
\label{eq:3-61}
\end{eqnarray}
where 
\begin{eqnarray*}
I'_{k}(\tau ) &:=& \left\{ \begin{array}{ccc}
\vspace{0.3cm}
\displaystyle \int_{0}^{\eta_{\delta}'}\left( 1-\tau \frac{-{\rm i}}{r-s_{0}\overline{z_{\alpha}}} \right)\frac{r^{\frac{2k-1}{2}}}{(r-s_{0}\overline{z_{\alpha}})^{2}}\exp\left( \frac{{\rm i}\tau}{r-s_{0}\overline{z_{\alpha}}} \right)\; {\rm d}r & {\rm if} & j\;\; {\rm is}\;\;{\rm odd},\\
\displaystyle \int_{0}^{\eta_{\delta}'}\left( 1-\tau \frac{{\rm i}}{r-s_{0}z_{\alpha}} \right)\frac{r^{\frac{2k-1}{2}}}{(r-s_{0}z_{\alpha})^{2}}\exp\left( \frac{-{\rm i}\tau}{r-s_{0}z_{\alpha}} \right)\; {\rm d}r & {\rm if} & j\;\; {\rm is}\;\;{\rm even}.
\end{array} \right.
\end{eqnarray*}

Now we provide a lemma for asymptotic behavior of $ I_{k}(\tau) $ as $\tau\longrightarrow\infty$;

\begin{lemma}[Lemma 3.5 in \cite{HIIS}, Lemma 3 in \cite{IIS}]\label{lem1-3}
Let $ n=1, 2, 3, \cdots $.
It holds for all $\alpha \in\, \bigr] -\frac{\pi}{2},\frac{\pi}{2} \bigl] $
\begin{eqnarray*}
\lefteqn{\hspace{1cm}
\lim_{\tau\longrightarrow\infty}
\tau^{\frac{2n+1}{2}}e^{-\frac{\tau}{2s_{0}}}e^{(-1)^{j}\frac{{\rm i}\tau\cos{\alpha}}{2s_{0}(1+\sin{\alpha})}}I_{n}(\tau)}\\
&&=(-1)^{j}{\rm i}s_{0}^{2n-1}2^{\frac{2n-1}{2}}(1+\sin{\alpha})^{\frac{2n-1}{2}}e^{(-1)^{j+1}{\rm i}\left(\frac{2n-1}{2}\right)\alpha}\Gamma\left(\frac{2n+1}{2}\right).
\end{eqnarray*}
\end{lemma}

\noindent This lemma implies that as $\tau\longrightarrow\infty$
\begin{eqnarray*}
\lefteqn{\hspace{1cm}
e^{-\frac{\tau}{2s_{0}}}I_{n}(\tau)}\\
&& \sim (-1)^{j}{\rm i}s_{0}^{2n-1}2^{\frac{2n-1}{2}}(1+\sin{\alpha})^{\frac{2n-1}{2}}e^{(-1)^{j+1}{\rm i}\left(\frac{2n-1}{2}\right)\alpha}e^{\frac{(-1)^{j+1}{\rm i}\tau\cos{\alpha}}{2s_{0}(1+\sin{\alpha})}}\Gamma\left(\frac{2n+1}{2} \right) \tau^{-\frac{2n+1}{2}}.
\end{eqnarray*}

\noindent Combining \eqref{eq:3-5} with \eqref{eq:3-6} and applying Lemma \ref{lem1-3} lead to \eqref{eq:2-9}.

As regards $I'_{n}(\tau)$, a similar method used in the proof of Lemma \ref{lem1-3} provided in \cite{HIIS, IIS} is applicable.
In consequence, as $\tau \longrightarrow 0 $, it holds for for all $\alpha \in\, \bigr] -\frac{\pi}{2},\frac{\pi}{2} \bigl] $
\begin{eqnarray*}
\lefteqn{e^{-\frac{\tau}{2s_{0}}}I'_{n}(\tau)\sim\left( 1+\tau\left(\frac{1}{2s_{0}}-(-1)^{j}\frac{\rm i}{2}\frac{\cos{\alpha}}{s_{0}(1+\sin{\alpha})}\right)\right)\times}\\
&& \times (-1)^{j}{\rm i}s_{0}^{2n-1}2^{\frac{2n-1}{2}}(1+\sin{\alpha})^{\frac{2n-1}{2}}e^{(-1)^{j+1}{\rm i}\left(\frac{2n-1}{2}\right)\alpha}e^{\frac{(-1)^{j+1}{\rm i}\tau\cos{\alpha}}{2s_{0}(1+\sin{\alpha})}}\Gamma\left(\frac{2n+1}{2} \right) \tau^{-\frac{2n+1}{2}}.
\end{eqnarray*}

\noindent From \eqref{eq:3-51} and \eqref{eq:3-61}, we obtain \eqref{eq:2-10}.

\subsection{Proof of Lemma \ref{lem1}}

The proof basically proceeds along the same line as that of Lemma 2 in \cite{II2007}, see also the proof of Lemma 3.4 in \cite{HIIS} in the case of Laplace equation.
However, due to multiple cracks in elasticity, we employ analytic continuation arguments on complex stress functions.

Now we assume that $A_{2n-1}^{(j)}=B_{2n-1}^{(j)}= 0$ for all $ n\geq 1 $.
Then, \eqref{eq:2-4} can be reduced to
\begin{eqnarray}
\mbox{\boldmath $u$}(r,\theta )=\sum_{k=0}^{\infty}\frac{A_{2k}^{(j)}}{2\mu}r^{k}\mbox{\boldmath $\varphi$}_{2k}(\theta )-\sum_{k=0}^{\infty}\frac{B_{2k}^{(j)}}{2\mu}r^{k}\mbox{\boldmath $\psi$}_{2k}(\theta ), \label{eq:3-14}
\end{eqnarray}
which means that $\mbox{\boldmath $u$}(r,\theta )$ is real analytic near the crack tip $(c_{j},c)$.
Moreover, it is easy to see that $\mbox{\boldmath $\sigma$}(\mbox{\boldmath $u$})\mbox{\boldmath $e$}_{2}=\mbox{\boldmath $0$}$ on $\{ r\in ]0, \epsilon_{0} [, \theta =0\}$ for a sufficiently small $\epsilon_{0}>0$.

Next, we construct the Goursat-Kolosov-Muskhelishvili stress functions $\phi_{\pm}(z)$ and $\omega_{\pm}(z)$, see \cite{M1977, IVT2011}, in each $B_{\pm}$, respectively.
Here we use notations $z=z_{1}+{\rm i}z_{2}:=x_{1}+{\rm i}(x_{2}-c)$,
\begin{eqnarray*}
B_{+} &:=& \left\{ z\; |\; z_{1}\in\, ]0, a[,\; z_{2}\in\, ]0, {\rm min}\{c, b-c \}[ \right\},\\
B_{-} &:=& \left\{ z\; |\; z_{1}\in\, ]0, a[,\; z_{2}\in\, ]-{\rm min}\{c, b-c \}, 0[ \right\}.
\end{eqnarray*}
Then $\phi_{\pm}(z)$ and $\omega_{\pm}(z)$ are holomorphic functions in $B_{\pm}$ of the complex variable $z$ due to the interior and boundary regularity results of the solution of $(*)$.
Moreover, it follows from the generalized Poincar\'e lemma (e.g. \cite{K2005}) and an argument in \cite{II2007, II2008, II2009} that $\phi_{\pm}(z), \omega_{\pm}(z)\in H^{1}(B_{\pm})$, respectively.
The displacement and the stress fields are given by the stress functions as follows
\begin{eqnarray*}
2\mu (u_{1} +{\rm i}u_{2}) &=& \kappa \phi_{\pm}(z)-\overline{\omega_{\pm}(z)}+(\overline{z} -z)\overline{\dot{\phi}_{\pm}(z)},\\
\sigma_{22}- {\rm i}\sigma_{12} &=& \phi'_{\pm}(z)+\overline{\dot{\omega}_{\pm}(z)}+(z-\overline{z})\overline{\ddot{\phi}_{\pm}(z)},
\end{eqnarray*}
where $\dot{\phi}(z) := d\phi /dz$.
The condition on $ \Sigma $, $ \sigma_{12}=\sigma_{22}=0 $, implies
\begin{eqnarray*}
\dot{\phi}_{+}(z_{1})+\overline{\dot{\omega}_{+}(z_{1})}=0, \quad \dot{\phi}_{-}(z_{1})+\overline{\dot{\omega}_{-}(z_{1})}=0 \quad {\rm on}\quad \Sigma .
\end{eqnarray*}
From this we define sectionally holomorphic functions $ \Psi_{1} (z) $ and $ \Psi_{2} (z) $ cut along $[0,a]\times\{c\}\setminus\Sigma$ in this way
\[
\Psi_{1} (z):=
\left\{ \begin{array}{ccc}
\vspace{0.3cm}
\dot{\phi}_{+}(z) & {\rm in} & B_{+}, \\
-\overline{\dot{\omega}_{+}(\overline{z})} & {\rm in} & B_{-},
\end{array}\right.
\quad
\Psi_{2} (z):=
\left\{ \begin{array}{ccc}
\vspace{0.3cm}
-\overline{\dot{\omega}_{-}(\overline{z})}& {\rm in} & B_{+}, \\
\dot{\phi}_{-}(z) & {\rm in} & B_{-}.
\end{array}\right.
\]
Meanwhile, $ \mbox{\boldmath $u$}^{+}=\mbox{\boldmath $u$}^{-} $ holds on $[0,a]\times\{c\}\setminus\Sigma$ and then we have
\begin{eqnarray*}
\frac{\kappa}{\mu}\dot{\phi}_{+}(z_{1})-\frac{1}{\mu}\overline{\dot{\omega}_{+}(z_{1})}=\frac{\kappa}{\mu}\dot{\phi}_{-}(z_{1})-\frac{1}{\mu}\overline{\dot{\omega}_{-}(z_{1})}\quad {\rm on}\quad [0,a]\times\{c\}\setminus\Sigma .
\end{eqnarray*}
Since both $ \sigma_{12}=\sigma_{22}=0 $ and $\mbox{\boldmath $u$}^{+}=\mbox{\boldmath $u$}^{-} $ holds on $ ]c_{j}, c_{j}+\epsilon_{0}[\times \{ c\} $ for a odd $j$ or $ ]c_{j}-\epsilon_{0}, c_{j}[\times \{ c\} $ for an even $j$, one sees
\[
\dot{\phi}_{+}(z_{1})=\dot{\phi}_{-}(z_{1}),\quad
\dot{\omega}_{+}(z_{1})=\dot{\omega}_{-}(z_{1}),\quad
\dot{\phi}_{\pm}(z_{1})=-\overline{\dot{\omega}_{\mp}(z_{1})}.
\]
By virtue of analytic continuation, we can define holomorphic functions $B_{+}\cup B_{-}$
\[
\dot{\phi} (z):=
\left\{ \begin{array}{ccc}
\vspace{0.3cm}
\dot{\phi}_{+}(z) & {\rm in} & B_{+}, \\
\dot{\phi}_{-}(z) & {\rm in} & B_{-},
\end{array}\right.
\quad
\dot{\omega} (z):=
\left\{ \begin{array}{ccc}
\vspace{0.3cm}
\dot{\omega}_{+}(z) & {\rm in} & B_{+}, \\
\dot{\omega}_{-}(z) & {\rm in} & B_{-}.
\end{array}\right.
\]
In the consequence it yields $\dot{\phi}(z)=-\overline{\dot{\omega}(\overline{z})}$, that is, $\Psi_{1}=\Psi_{2}$.
Then we sees $\mbox{\boldmath $\sigma$}(\mbox{\boldmath $u$})\mbox{\boldmath $e$}_{2}=\mbox{\boldmath $0$} $ on $]0,a[\times \{c\}$.
Moreover, from the assumption {\rm ($\dag$)} for $\mbox{\boldmath $g$}$, $\mbox{\boldmath $u$}$ has at least $H^2$ regularity near all the corner points $ \mbox{\boldmath $O$} $, $ (a,0) $, $ (0,c) $, $(a,c)$, $(0,b)$ and $(a,b)$, e.g. \cite{KMR, II2009}.
Therefore, one sees that the restriction of $\mbox{\boldmath $u$} $ to $\Omega^{-}$ is in $ H^{2}(\Omega^{-})$ and satisfies
\[
\left\{ \begin{array}{l}
\vspace{0.3cm}
\mu\triangle\mbox{\boldmath $u$} + (\lambda +\mu) \nabla (\nabla\cdot\mbox{\boldmath $u$})=\mbox{\boldmath $0$} \quad {\rm in} \quad \Omega^{-},\\
\vspace{0.3cm}
\mbox{\boldmath $\sigma$}(\mbox{\boldmath $u$})\mbox{\boldmath $\nu$}=\mbox{\boldmath $0$} \quad {\rm on} \quad ]0,a[\,\times \{c\},\\
\mbox{\boldmath $\sigma$}(\mbox{\boldmath $u$})\mbox{\boldmath $\nu$}=\mbox{\boldmath $g$} \quad {\rm on} \quad \partial\Omega\cap \{x_{2}<c\}.
\end{array}\right.
\]
Then, the divergence theorem yields that for an arbitrary $\mbox{\boldmath $\rho$}\in {\cal R} $
\begin{eqnarray*}
0 = \int_{\Omega^{-}}\mbox{\boldmath $\rho$}\cdot (\mu\triangle\mbox{\boldmath $u$} + (\lambda +\mu) \nabla (\nabla\cdot\mbox{\boldmath $u$}))\; {\rm d}\mbox{\boldmath $x$} = \int_{\partial\Omega\cap \{x_{2}<c\}}\mbox{\boldmath $\rho$}\cdot\mbox{\boldmath $g$}\; {\rm d}s_{\mbox{\boldmath $x$}}.
\end{eqnarray*}
Since $ \mbox{\boldmath $g$} $ satisfies the condition {\rm ($\dag$)}, it is a contradiction. 

Similarly, the above argument is valid for the case of $\Omega^{+}$.
The proof of Lemma \ref{lem1} is completed.

\section{Conclusion}
In this paper, an inverse crack problem in a two dimensional linearized elasticity is considered.
Applying the enclosure method combined with the Kelvin transform, we established two types of extraction formulae \eqref{eq:1-10} and \eqref{eq:1-11} which enable us to detect multiple linear cracks located on a line by processing mathematically a single set of observation data. 
Formula \eqref{eq:1-10} including (S1) in Theorem \ref{th0} is commonly used and is extended the result of the case of electric conductive body \cite{HIIS, IIS}.
Although outline of the proof of \eqref{eq:1-10} is given in \cite{HI2020}, in this paper we provided the proof of key lemma to show \eqref{eq:1-10}.
Formula \eqref{eq:1-11} is derived by applying the idea of taking logarithmic derivative of the indicator function \cite{I2011}, which makes us possible to obtain more information of unknown cracks.
However, these are only theoretical results, therefore it is important to demonstrate the feasibility of the method.
Computational implementation of the method is expected as our next problem.

\section*{Acknowledgments}
MI was partially supported by a Grant-in-Aid for Scientific Research (C) (No. 17K05331) and (B)(No. 18H01126) of the Japan Society for the Promotion of Science (JSPS).
HI was partially supported by a Grant-in-Aid for Scientific Research (C) (No. 18K03380) of JSPS.
This work is supported by JSPS and the Russian Foundation for Basic Research (RFBR) under the Japan - Russia Research Cooperative Program (project No. JPJSBP120194824).

\end{document}